\documentclass[twoside,11pt]{article}
\title{} \author{} \date{}
\usepackage{latexsym,amssymb,times}
\newtheorem{te}{Theorem}[section]

\newtheorem{fac}[te]{Fact}

\newtheorem{cla}[te]{Claim}

\def\dok{\noindent{\bf Proof. }}
\def\kdok{\hfill $\Box$ \par \vspace*{2mm} }
\def\a{\alpha}
\def\b{\beta}
\def\g{\gamma}
\def\o{\omega}

\def\f{\varphi}

\def\l{\lambda}

\def\r{\rho}
\def\s{\sigma}
\def\t{\tau}

\def\A{\mathbb A}
\def\B{\mathbb B}
\def\BD{\mathbb D}
\def\BI{\mathbb I}
\def\G{\mathbb G}
\def\N{\mathbb N}
\def\P{\mathbb P}
\def\R{\mathbb R}
\def\X{\mathbb X}
\def\Y{\mathbb Y}
\def\S{\mathbb S}
\def\Q{\mathbb Q}
\def\BT{\mathbb T}

\def\Z{\mathbb Z}

\def\CP{{\mathcal P}}
\def\T{{\mathcal T}}
\def\la{\langle}
\def\ra{\rangle}

\def\Aut{\mathop{\rm Aut}\nolimits}

\def\sq{\mathop{\rm sq}\nolimits}
\def\ro{\mathop{\rm ro}\nolimits}
\def\Fin{\mathop{\rm Fin}\nolimits}

\def\ar{\mathop{\rm ar}\nolimits}
\def\sh{\mathop{\rm sh}\nolimits}
\def\Rado{\mathop{\G _{{\mathrm{Rado}}}}\nolimits}
\def\Th{\mathop{\rm Th}\nolimits}
\begin{document}
\thispagestyle{plain}
\begin{center}
           {\large \bf \uppercase{Posets of copies of countable ultrahomogeneous tournaments}}
\end{center}
\begin{center}
{\bf Milo\v s S.\ Kurili\'c\footnote{Department of Mathematics and Informatics, Faculty of Science, University of Novi Sad,
              Trg Dositeja Obradovi\'ca 4, 21000 Novi Sad, Serbia.
              email: milos@dmi.uns.ac.rs}
and Stevo Todor\v cevi\'c\footnote{Department of Mathematics, University of Toronto, Canada;
              Institut de Math\'{e}matiques de Jussieu, CNRS, Paris, France; and
              Matemati\v cki Institut, SANU, Belgrade, Serbia.
              e-mail: stevo@math.toronto.edu, stevo.todorcevic@imj-prg.fr, stevo.todorcevic@sanu.ac.rs}}
\end{center}
\begin{abstract}
\noindent
The {\it poset of copies} of a relational structure ${\mathbb X}$
is the partial order ${\mathbb P} ({\mathbb X} ) := \langle \{ Y \subset X: {\mathbb Y} \cong {\mathbb X}\}, \subset \rangle$
and each similarity of such posets (e.g.\ isomorphism, forcing equivalence = isomorphism of Boolean completions, ${\mathbb B} _{\mathbb X}:={\mathrm{ro}}\,{\mathrm{sq}}\, {\mathbb P} ({\mathbb X} )$)
determines a classification of structures.
Here we consider the structures from Lachlan's list of countable ultrahomogeneous tournaments:
${\mathbb Q} $ (the rational line), ${\mathbb S} (2)$ (the circular tournament),
and ${\mathbb T} ^\infty$ (the countable homogeneous universal tournament);
as well as the ultrahomogeneous digraphs ${\mathbb S} (3)$, ${\mathbb Q} [{\mathbb I}_n]$, ${\mathbb S} (2)[{\mathbb I}_n]$ and ${\mathbb T} ^\infty [{\mathbb I}_n]$
from Cherlin's list.

If $\G _{{\mathrm{Rado}}}$ (resp.\ ${\mathbb Q} _n$) denotes the countable homogeneous universal graph (resp.\ $n$-labeled linear order),
it turns out that ${\mathbb P} ({\mathbb T} ^\infty)\cong {\mathbb P} (\G _{{\mathrm{Rado}}})$
and that ${\mathbb P} ({\mathbb Q} _n)$ densely embeds in ${\mathbb P} ({\mathbb S} (n))$, for $n\in\{ 2,3\}$.

Consequently, ${\mathbb B} _{\mathbb X} \cong {\mathrm{ro}}\, ({\mathbb S} \ast \pi)$,
where ${\mathbb S}$ is the poset of perfect subsets of $\R$
and $\pi$ an ${\mathbb S}$-name such that $1_{\mathbb S} \Vdash `` \pi $ is a separative, atomless and $\sigma$-closed forcing"
(thus $1_{\mathbb S} \Vdash `` \pi \equiv_{forc}(P(\omega )/{\mathrm{Fin}} )^+$", under CH),
whenever ${\mathbb X}$ is a countable structure equimorphic with ${\mathbb Q}$, ${\mathbb Q} _n$, ${\mathbb S} (2)$, ${\mathbb S} (3)$, ${\mathbb Q} [{\mathbb I}_n]$ or ${\mathbb S} (2)[{\mathbb I}_n]$.

Also, ${\mathbb B} _{\mathbb X} \cong {\mathrm{ro}}\, ({\mathbb S} \ast \pi)$,
where  $1_{\mathbb S} \Vdash ``\pi $ is an $\o$-distributive forcing",
whenever ${\mathbb X}$ is a countable graph embedding $\G _{{\mathrm{Rado}}}$,
or a countable tournament embedding ${\mathbb T} ^\infty$,
or ${\mathbb X} ={\mathbb T} ^\infty [{\mathbb I}_n]$. 
\end{abstract}
\section{Introduction}\label{S1}
If $\X $ is a relational structure and $\P (\X ) = \{ Y \subset X: \Y \cong \X\}$ the set of copies of $\X$ inside $\X$,
the partial ordering $\la \P (\X ), \subset \ra$ will be called the {\it poset of copies of $\X$} and shortly denoted by
$\P (\X )$, whenever the context admits.

It is easy to see that the correspondence $\X \mapsto \B _\X$ (where $\B _\X$ is the Boolean completion
of the separative quotient of the poset $\P (\X )$, $\ro\sq \P (\X )$) extends to a functor from the category of all
relational structures and isomorphisms  to its subcategory of all homogeneous complete Boolean algebras and,
defining two relational structures $\X $ and $\Y$ to be similar iff $\B _\X\cong \B _\Y$, we obtain a coarse classification of relational structures
(see \cite{KurS}). The position of this similarity in the hierarchy of set-theoretical and model-theoretical similarities of structures
was investigated in \cite{Kdif,KuMo};
in particular, for relational structures $\X$ and $\Y$ we have:
\begin{equation}\label{EQ2372}
\X \rightleftarrows \Y \Rightarrow \P (\X )\equiv _{forc}\P (\Y ) \Leftrightarrow \B _\X \cong \B _\Y ,
\end{equation}
where $\rightleftarrows$ denotes the equimorphism (bi-embedability) relation and $\equiv _{forc}$ the forcing equivalence of posets.
So, the mentioned classification of structures can be explored
using the methods of set-theoretic forcing.

In this paper we continue the investigation of countable ultrahomogeneous relational structures in this context.
By (\ref{EQ2372}), a statement concerning a countable ultrahomogeneous structure $\X$ holds for all the structures from its equimorphism class.
For example, if $\Q$ denotes the rational line, $\la Q ,<_\Q\ra$, then $\B _\Q\cong \B _\X$, for each
countable non-scattered linear order $\X$.

All the definitions and facts concerning ultrahomogeneous structures used in this paper can be found in the survey \cite{Mac} of Macpherson.
By $\Rado$ we denote the Rado graph and by $\Q _n$ (for $n\in \N$) the countable ultrahomogeneous $n$-labeled linear order, that is the structure
$\Q _{n} :=\la Q, < _\Q ,A_1 ,\dots ,A_n\ra$, where $\{ A_1 ,\dots ,A_n\} $ is a partition of the set $Q$ such that the sets $A_i$, $i\leq n$, are dense in $\Q$.

In order to state the known results which will be used in this paper, by $\S$ we denote the Sacks perfect set forcing
(the set of perfect subsets of $\R$ ordered by the inclusion)
and, in order to avoid repetition, we introduce the following notation for two properties of a countable relational structure $\X$:
\begin{itemize}
\item[$\CP _1$:] $\P (\X ) \equiv_{forc} \S \ast \pi$,
where $\pi$ is an $\S$-name for a preorder  and $1_\S \Vdash `` \pi $ is a separative, atomless and $\sigma$-closed forcing";
\item[$\CP _2$:] $\P (\X ) \equiv_{forc} \S \ast \pi$,
where  $\pi$ is an $\S$-name for a preorder and $1_\S \Vdash ``\pi $  is an $\o$-distributive forcing".
\end{itemize}
\begin{fac}\label{T2659}
Let $\sh (\S )$ denote the size of the continuum in the Sacks extension (the cardinal
$\kappa$ such that $1_\S \Vdash {\mathfrak c}=\check{\kappa }$) and let $\X$ be a countable relational structure.

(a) $\CP _1$ implies $\CP _2$;

(b) If $\CP _1$ is true and $\sh (\S )=\aleph _1$, then $1_\S \Vdash `` \pi \equiv_{forc}(P(\o )/\Fin )^+$";

(c) CH and, more generally, the equality ${\mathfrak b}=\aleph _1$ implies that $\sh (\S )=\aleph _1$.
\end{fac}
\dok
Since each $\s$-closed forcing is $\o$-distributive (a) is true.
It is a folklore fact that under CH each separative, atomless and $\sigma$-closed forcing of size ${\mathfrak c}$
is forcing equivalent to $(P(\o )/\Fin )^+$, which proves (b). For (c) see \cite{Sim}.
\hfill $\Box$
\begin{te} \label{T4120}
(a) Each countable linear order embedding $\Q$ has property $\CP _1$ \cite{KurTod}.

(b) Each countable $n$-labeled linear order embedding $\Q_n$ has property $\CP _1$ \cite{KurTod3}.

(c) Each countable graph embedding $\Rado$ has property $\CP _2$ \cite{GuzTod,KurTod2,KurTod1}.\footnote{
In \cite{KurTod2} and \cite{KurTod1}it was proved that
$\P (\Rado ) \equiv_{forc} \P \ast \pi$,
where  $\P$ is a poset which adds a generic real,
has the 2-localization property (and, hence, the Sacks property)
has the $\aleph _0$-covering property (thus preserves $\o _1$)
and does not produce splitting reals
and $\pi$ is a $\P$-name for a preorder such that $1_\P \Vdash ``\pi \mbox{ is an }\o\mbox{-distributive forcing}"$.
The forcing equivalence $\P (\Rado ) \equiv_{forc} \S \ast \pi$ from $\CP _2$ was proved in \cite{GuzTod}.
}
\end{te}
The aim of this paper is to complete the picture for countable ultrahomogeneous tournaments.
We recall Lachlan's classification of these structures \cite{Lach}:
Each countable ultrahomogeneous tournament is isomorphic to one of the following:
\begin{itemize}
\item ${\mathbb Q} $, the rational line,
\item  $\S (2)$, the dense local order (the circular tournament),
\item  ${\mathbb T} ^\infty$, the countable random (i. e.\ homogeneous universal) tournament.
\end{itemize}
\noindent
In Sections \ref{S2} and \ref{S3} we show that ${\mathbb T} ^\infty$ has $\CP _2$ and that $\S (2)$ has $\CP _1$ and
in Section \ref{S4} we obtain similar results for infinitely many ultrahomogeneous digraphs
from Cherlin's list \cite{Che2}: $\S (3)$, ${\mathbb T} [I_n]$ and $I_n[{\mathbb T} ]$, where ${\mathbb T} \in \{ \Q , {\mathbb T} ^\infty, \S (2)\}$ and $n\in \N$.
More precisely, the main results of the paper are the following.
\begin{itemize}
\item $\P (\BT ^\infty)\cong \P (\Rado )$ and, hence, $\B _{\BT ^\infty}\cong \B _{\Rado }$.\\
      Each countable tournament $\X$ embedding $\BT ^\infty$ has property $\CP _2$.
\item $\P (\Q _2)$  densely embeds in $\P (\S (2))$ and, hence, $\B _{\S (2)}\cong \B _{\Q _2}$.\\
      Each countable tournament $\X$ equimorphic with $\S (2)$  has property $\CP _1$.
\item $\P (\Q _3)$  densely embeds in $\P (\S (3))$ and, hence, $\B _{\S (3)}\cong \B _{\Q _3}$. \\
      Each countable digraph $\X$ equimorphic with $\S (3)$ has property $\CP _1$.
\end{itemize}
The following elementary fact will be used in the sequel.
\begin{fac}\label{T2658}
Let $\X=\la X,\r \ra$ be a countable ultrahomogeneous relational structure of a finite language. Then

(a) The theory $\Th (\X )$ is $\o$-categorical and admits quantifier elimination;

(b) $\P (\X )$ is equal to the set of domains of elementary substructures of $\,\X$.
\end{fac}
\dok
For (a) see \cite{Hodg}, p.\ 350. If $A\in \P (\X )$, then $\A \models \Th (\X )$ and, since by (a)
$\Th (\X )$ is model complete, $\A \prec \X$. Conversely, if $\A =\la A,\r\upharpoonright A \ra \prec \X$, then $\A \equiv X$ and, since
$\Th (\X )$ is $\o$-categorical, $\A \cong \X$, that is, $A\in \P (\X )$.
\hfill $\Box$
\section{The random tournament}\label{S2}
\paragraph{The Rado graph}
If $\la G, \sim \ra$ is a graph  and $K\subset H \in [G]^{<\o }$, let us define
$$
G^H_K \!:=\Big\{ v\in G \setminus H : \forall k\in K \, (v\sim k) \; \land \; \forall h\in H\setminus K \, (v\not\sim h)\Big\}.
$$
(Clearly, $G^\emptyset _\emptyset =G$.)
The {\it Rado graph}, $\Rado$, \cite{Rado} (the Erd\H{o}s-R\'enyi graph \cite{Erdos2}, the countable random graph) is
the unique (up to isomorphism) countable homogeneous universal graph
and the Fra\"{\i}ss\'e limit of the amalgamation class of all finite graphs; see \cite{Fra}, where a proof of the following fact
can be found.
\begin{fac}\label{T2645}
For a countable graph $\G =\langle G, \sim \rangle $ the following is equivalent

{\rm{(g1)}} $\G \cong \Rado$,

{\rm{(g2)}} $G^H_K\neq \emptyset$, whenever $K\subset H\in [G]^{<\o }$,

{\rm{(g3)}} $|G^H_K|=\o$, whenever $K\subset H\in [G]^{<\o }$.
\end{fac}
\paragraph{The random tournament}
If $\la T, \rightarrow \ra$ is a tournament, and $K\subset H \in [T]^{<\o }$, let
$$T^H_K \!:=\Big\{ v\in T \setminus H : \forall k\in K \; (k\rightarrow v) \; \land \;
\forall h\in H\setminus K \; (v\rightarrow h)\Big\}.
$$
(Clearly, $T^\emptyset _\emptyset =T$.)
The {\it random tournament}, $\BT ^\infty$, is the unique (up to isomorphism) countable homogeneous universal tournament
and the Fra\"{\i}ss\'e limit of the amalgamation class of all finite tournaments (see \cite{Fra}).
\begin{fac}\label{T2641}
For a countable  tournament $\BT =\la T, \rightarrow \ra$ the following is equivalent

{\rm{(t1)}} $\BT \cong \BT ^\infty$,

{\rm{(t2)}} $T^H_K\neq \emptyset$, whenever $K\subset H\in [T]^{<\o }$,

{\rm{(t3)}} $|T^H_K|=\o$, whenever $K\subset H\in [T]^{<\o }$.
\end{fac}
\dok
(t1) $\Rightarrow$ (t2). Let $\BT =\la T, \rightarrow \ra\cong \BT ^\infty$, $K\subset H \in [T]^{<\o }$ and $p\not\in T$.
Then $\BT _0 :=\la  H\cup \{ p \} , \r \ra$, where
$$
\r = (\rightarrow \upharpoonright H)
\cup \{ \la k,p \ra : k\in K \}
\cup \{ \la p,h \ra : h\in H\setminus K \} ,
$$
is a finite tournament
and, since the age of $\BT$ is the class of all finite tournaments, there is an embedding $f :\BT _0 \hookrightarrow \BT$.
Now the restriction $\f := f^{-1}\upharpoonright f[H] $ is a finite partial isomorphism of $\BT$ which maps $f[H]$ onto $H$ and, by the ultrahomogeneity of $\BT$ there is $F\in \Aut (\BT)$ such that $\f \subset F$. Let $v:= F(f(p))$. For $k\in K$ we have
$\la k,p \ra \in \r$ and, hence,  $\la f(k),f(p) \ra \in \rightarrow$, which implies $\la F(f(k)),F(f(p)) \ra \in \rightarrow$.
Since $F(f(k))=\f (f(k))= f^{-1}(f(k))=k$, we have $\la k,v \ra\in \rightarrow$. Similarly, $\la v,h \ra\in \rightarrow$,
for all $h\in H\setminus K$, and, thus, $v\in T^H_K$.

(t2) $\Rightarrow$ (t3). Suppose that (t2) is true and that $T_K^H=\{ v_1,\dots ,v_n\}$.
Then, by (t2) there is $v\in T_K^{H\cup \{ v_1,\dots ,v_n\} }$ and, hence, $v\in T_K^H$ and $v\not\in H\cup \{ v_1,\dots ,v_n\}$,
which is a contradiction.

(t3) $\Rightarrow$ (t1). Assuming (t3) we show first that for each $n\in \N $ each finite tournament $\A$ of size $n$ embeds in $\BT$. For $n=1$ the statement is obviously true. Suppose that it is true for $n$ and that
$\A=\la A, \r\ra$ is a tournament, where $A=\{ a_1 , \dots a_{n+1} \}$. Then for $\A _0=\la A_0, \r \upharpoonright A_0\ra$,
where $A_0=\{ a_1 , \dots a_{n} \}$, there is an embedding $f: \A _0 \hookrightarrow \BT$ and if
$K:=\{  i\leq n : \la a_i , a_{n+1}\ra \in \r \}$, by  (t3) there is $v\in T$ such that
$f(a_i)\rightarrow v$, for each $i\in K$, and $v\rightarrow f(a_i)$, for all $i\in \{ 1, \dots , n \}\setminus K$.
Thus $f[A_0]\cup \{ v \}$ is a copy of $\A$
in $\BT$.

Now we show that $\BT$ has the 1-extension property. Let $\f :H\rightarrow T$ be a finite partial isomorphism,
$v\in T\setminus H$ and $K:=\{ k\in H : k\rightarrow v\}$. By (t3) there is $w\in T$ such that
$\f (k)\rightarrow w$, for all $k\in K$ and $w\rightarrow \f (h)$, for all $h\in H\setminus K$.
Thus $\f \cup \{ \la v,w \ra \}$ is a finite partial isomorphism of $\BT$.
\hfill $\Box$
\begin{te}\label{T2660}
$\P (\BT ^\infty)\cong \P (\Rado )$ and, hence, $\B _{\BT ^\infty}\cong \B _{\Rado }$.

Each countable tournament embedding $\BT ^\infty$ has property $\CP _2$.
\end{te}
\dok
W.l.o.g\ we suppose that $\Rado =\la \o , \sim \ra$ and define a binary relation $\rightarrow$ on the set $\o$ in the following way:
for $m,n\in \o$ let
\begin{equation}\label{EQ2369}
m \rightarrow n \Leftrightarrow ( m<n \land m\sim n) \lor (m>n \land m\not\sim n).
\end{equation}
Since the relations $\sim$ and $\sim ^c := \o ^2 \setminus \sim $ are symmetric, by (\ref{EQ2369})   we have
\begin{equation}\label{EQ2370}
\rightarrow = ( < \cap \sim ) \cup (<^{-1} \cap \sim ^c ) \;\;\mbox{ and }\;\;
\rightarrow^{-1} = ( <^{-1}\cap \sim  ) \cup ( < \cap  \sim ^c)
\end{equation}
Now we have: $\rightarrow \cap \Delta _\o =\emptyset$ so the relation $\rightarrow$ is irreflexive,
$\rightarrow \cap \rightarrow^{-1}=\emptyset$, and $\rightarrow$ is asymmetric and
$\rightarrow \cup \rightarrow^{-1}=< \cup <^{-1}=\o ^2\setminus \Delta _\o$; thus
the structure $\BT :=\la \o ,\rightarrow\ra$ is a tournament.

For a proof that $\BT \cong \BT ^\infty$, we check (t2) of Fact \ref{T2641}. If $K\subset H \in [\o ]^{<\o }$,
then, by (g3) of Fact \ref{T2645}, there is $v\in G_K^H$ such that $v> h$, for all $h\in H$.
Now, if $k\in K$, then $k<v$ and $k\sim v$ so, by (\ref{EQ2369}), $k\rightarrow v$.
If $h\in H\setminus K$, then $v>h$ and $v\not\sim h$ so, by (\ref{EQ2369}) again, $v\rightarrow h$.
Thus $v\in T_K^H \neq \emptyset$,  (t2) is true and $\BT \cong \BT ^\infty$ indeed.

In order to prove that
\begin{equation}\label{EQ2371}
\P (\la \o , \sim \ra) =\P (\la \o ,\rightarrow\ra)
\end{equation}
we take first $A\in \P (\la \o , \sim \ra)$ and show that the countable tournament $\la A ,\rightarrow \upharpoonright A\ra$
satisfies (t2) of Fact \ref{T2641}. So, if $K\subset H \in [A ]^{<\o }$, then, since $\la A ,\sim \upharpoonright A\ra\cong \Rado$,
by (g3) of Fact \ref{T2645}, there is $v\in A_K^H$ such that $v> h$, for all $h\in H$.
Now, if $k\in K$, then $k<v$ and $k\sim v$ so, by (\ref{EQ2369}), $k\rightarrow v$.
If $h\in H\setminus K$, then $v>h$ and $v\not\sim h$ so, by (\ref{EQ2369}) again, $v\rightarrow h$.
Thus  $\la A ,\rightarrow \upharpoonright A\ra$ satisfies (t2) and $\la A ,\rightarrow \upharpoonright A\ra \cong \BT ^\infty$, which means that
$A\in \P (\la \o ,\rightarrow\ra)$.

Conversely, we take $A\in \P (\la \o , \rightarrow \ra)$ and show that the graph $\la A ,\sim \upharpoonright A\ra$
satisfies (g2) of Fact \ref{T2645}. So, if $K\subset H \in [A ]^{<\o }$, then, since $\la A ,\rightarrow \upharpoonright A\ra\cong \BT ^\infty$,
by (t3) of Fact \ref{T2641}, there is $v\in A_K^H$ such that $v> h$, for all $h\in H$.
If $k\in K$, then $k<v$ and $k\rightarrow v$ so, by (\ref{EQ2369}), $k\sim v$, that is $v\sim k$.
If $h\in H\setminus K$, then $v>h$ and $v\rightarrow h$ so, by (\ref{EQ2369}) again, $v\not\sim h$.
Thus  $\la A ,\sim \upharpoonright A\ra$ satisfies (g2) and $\la A ,\sim \upharpoonright A\ra \cong \Rado$,
which means that
$A\in \P (\la \o ,\sim\ra)$.

So, since $\la \o ,\rightarrow\ra\cong\BT ^\infty$ by (\ref{EQ2371}) we have $\P (\BT ^\infty)\cong \P (\Rado )$
and, hence, $\P (\BT ^\infty)\equiv_{forc} \P (\Rado )$.
If $\X$ is a countable tournament and $\BT ^\infty \hookrightarrow \X$, then, by the universality of $\BT ^\infty$, $\X \hookrightarrow\BT ^\infty$, so
$\X \rightleftarrows\BT ^\infty $ and, hence, $\P (\X )\equiv _{forc}\P (\BT ^\infty)\equiv _{forc} \P (\Rado )$  which, together with Theorem \ref{T4120}(c)
implies that $\X$ has property $\CP _2$.
\hfill $\Box$
\section{The dense local order}\label{S3}
\paragraph{The countable homogeneous universal $n$-labeled linear order}
For $n\in \N$ let $L_n=\la R, \a _1 , \dots, \a _n  \ra$ be a relational language, where
$\ar (R)=2$ and $\ar (\a _i)=1$, for $i\leq n$. We recall that the $L_n$-structures
of the form $\X =\la X, <, A_1 ,\dots ,A_n\ra$, where $<$ is a linear order on the set $X$ and $\{A_1 ,\dots ,A_n\}$
a partition of $X$, are called {\it $n$-labeled linear orders}.
Since the $L_n$-structure $\Q _{n}$
is ultrahomogeneous, the $L_n$-theory $\T _n$ saying that an $L_n$-structure
$\X =\la X, <, A_1 ,\dots ,A_n\ra$ is
a model of $\T _n$ iff $\la X, <\ra$ is a dense linear order without end-points and  $\{A_1 ,\dots ,A_n\}$
a partition of $X$ into dense subsets of $\la X, <\ra$ is $\o$-categorical. Consequently we have  $D\in \P (\Q _n)$ iff
$\la D,  <_\Q \upharpoonright D, A_1 \cap D , \dots , A_n \cap D\ra \models \T _n$, that is
\begin{fac}\label{T2653}
$D\in \P (\Q _n)$ if and only if  $\la D, <_\Q \upharpoonright D\ra$ is dense linear order without end points and the sets $A_i\cap D$, for $i\in \{ 1,\dots ,n\}$,
are its dense subsets.
\end{fac}
\paragraph{The dense local order $\S (2)$}
If $q_1,q_2 \in Q$ and $q_1\neq q_2$, then, since $q_1-q_2\neq k\pi$, for all $k\in \Z$, $e^{q_1 i}$ and $e^{q_2 i}$ are different and non-antipodal points of the unit circle $S^1:=\{ e^{ti}: t\in [0,2\pi )\}$ in the complex plane and $S= \{ e^{qi}: q\in Q \}$ is a dense subset of $S^1$.
The {\it dense local order} is the tournament $\S (2)=\la S, \rightarrow \ra$, where
\begin{equation}\label{EQ2357}\textstyle
e^{q_1 i}\rightarrow e^{q_2 i} \;\;\Leftrightarrow \;\; q_2 -q_1 \in \bigcup _{k\in \Z }(2k\pi ,2k\pi +\pi  ),
\end{equation}
which means that the shorter oriented path from $e^{q_1 i}$ to $e^{q_2 i}$ is the anticlockwise oriented one.
In order to simplify notation let $L_2=\la R , \a ,\b \ra$.

Clearly,  $\{ A , B\}$ is a partition of the set $S$ into the left and right part, where
\begin{eqnarray*}\textstyle
A & := &\textstyle\Big\{ e^{qi}: q\in \bigcup _{k\in \Z}(\frac{\pi}{2}+2k\pi , \frac{3\pi}{2}+2k\pi) \cap Q\Big\} \;\mbox{ and }\;\\
B & := &\textstyle\Big\{ e^{qi}: q\in \bigcup _{k\in \Z}(\frac{3\pi}{2}+2k\pi , \frac{5\pi}{2}+2k\pi)\cap Q\Big\}.
\end{eqnarray*}
So $\la S, \rightarrow , A ,B \ra$ is an $L_2$-structure
and the $L_2$-formula
\begin{eqnarray}
\l (u, v) & := & \Big[ \Big(( \a (u) \land \a (v) ) \lor ( \b (u) \land \b (v))\Big)\land R(u,v) \Big]  \nonumber \\
      &      \lor       & \Big[ \Big(( \a (u) \land \b (v))\lor (\b (u) \land \a (v))\Big)\land R(v,u) \Big]\label{EQ2356}
\end{eqnarray}
defines the tournament relation $\rho := \{ \la x,y \ra \in S^2: \la S, \rightarrow , A ,B \ra \models \l [x,y]\}$ on the set $S$,
which preserves $\rightarrow$ between the elements of the same part,
and reverses $\rightarrow$ between the elements of different parts, namely,
\begin{equation}\label{EQ2360}
\r = \Big[\rightarrow \cap \;\Big((A\times A)\cup (B\times B)\Big)\Big] \cup \Big[\rightarrow ^{-1} \cap \;\Big((A\times B)\cup (B\times A)\Big)\Big].
\end{equation}
It is easy to see that $\la S, \rho \ra$ is a dense linear order without end points
and that $A$ and $B$ are its dense subsets (see \cite{Lach}, p.\ 434),
which means that the $L_2$-structure $\Y := \la S, \rho , A ,B \ra$ is a model of $\T_2$
and, since $\T_2$ is an $\o$-categorical theory, $\Y\cong\Q _2$.

For $x,y\in S^1$, let $x^\smallfrown y$ denote the set of elements of $S$ belonging to the shorter arc determined by $x$ and $y$
and let $a(x)$ denote the antipodal point of $x$.
\begin{te}\label{T2373}
$\P (\Q _2)$  densely embeds in $\P (\S (2))$ and, hence, $\B _{\S (2)}\cong \B _{\Q _2}$.

Each countable tournament equimorphic with $\S (2)$  has property $\CP _1$.
\end{te}
\dok
Since $\Y := \la S, \rho , A ,B \ra\cong \Q _2$ we have $\P (\Q _2) \cong \P (\Y)$ and we show that $\P (\Y )$ is a dense subset of $\P (\S (2))$.
First we prove that $\P (\Y ) \subset \P (\S (2))$.
So, if $D\in \P (\Y )$,  then there is an isomorphism
\begin{equation}\label{EQ2366}
F: \la S, \rho , A ,B \ra \rightarrow_{iso} \la D, \rho \upharpoonright D , A\cap D ,B\cap D \ra
\end{equation}
and in order to prove that $D\in \P (\S (2))$ it remains to be shown that the mapping
$F: \la S, \rightarrow \ra \rightarrow \la D, \rightarrow \upharpoonright D \ra$ is an isomorphism.
Clearly, the relation $\rightarrow$ is defined by the formula $\l$ in the structure $\Y$, that is
\begin{equation}\label{EQ2367}
\forall x,y  \in S \;\; \Big(x \rightarrow y \Leftrightarrow \la S, \rho , A,B \ra \models \l [x,y] \Big).
\end{equation}
Now for $x,y \in S$ we have:
$x\rightarrow y$
iff (by (\ref{EQ2367})) $\la S, \rho , A ,B \ra  \models \l [x,y] $
iff (by (\ref{EQ2366})) $\la D, \rho \upharpoonright D , A\cap D ,B\cap D \ra  \models \l [F(x),F(y)] $
iff (since $\l$ is a $\Sigma _0$-formula and, thus $(D,S)$-absolute) $\la S, \rho , A ,B \ra  \models \l [F(x),F(y)] $
iff (by (\ref{EQ2367})) $F(x)\rightarrow F(y)$.
Thus $F: \la S, \rightarrow \ra \rightarrow \la D, \rightarrow \upharpoonright D \ra$ is an isomorphism, $D\in \P (\S(2))$ and we have proved that
$\P (\Y ) \subset \P (\S (2))$.
\begin{cla}\label{T2652}
If $D\in \P (\S (2))$, then $\la D ,\r \upharpoonright D\ra$ is a dense linear order with at most one end point and
$A_1:=A\cap D$ and $B_1:=B\cap D$ are its dense subsets.
\end{cla}
\dok
By Fact \ref{T2658}(b),
$D\in \P (\S (2))$ implies that $\BD:=\la D ,\rightarrow  \upharpoonright D\ra$ is an elementary substructure of $\S(2)$. So, by the Tarski-Vaught theorem,
in particular, for each $L_b$-formula $\theta (u,v,w)$ we have:
\begin{equation}\label{EQ2364}
\forall x,y \in D \;
\Big( \exists s\in S \;\; \S(2) \models  \theta [x,y,s]
\;\Rightarrow \exists z\in D \;\; \BD \models \theta [x,y,z]   \Big).
\end{equation}
Now $\la D,\r \upharpoonright D \ra$ is a linear order and we prove that $A_1$ is its dense subset, that is
\begin{equation}\label{EQ2365}
\forall x,y \in D \;\; \Big(x\r y \Rightarrow \exists z\in A_1 \;\; x\r  z\r y\Big).
\end{equation}
So, let $x,y \in D$ and $x\r y$. Then, since $\rho$ is a strict linear order, $\neg y \r x$.

If $x,y \in A_1$,
then by (\ref{EQ2360}) we have $x \rightarrow y$.
Since for $s\in x^\smallfrown y$ we have $x \rightarrow s \rightarrow y$,
by (\ref{EQ2364}) there is $z\in D$ such that $x \rightarrow z \rightarrow y$.
Now $z\in B_1$ would imply that $y \r z \r x$ and, hence, $y \r x$, which is false.
Thus $z\in A_1$
and, by (\ref{EQ2360}), $x \r z \r y$.

If $x,y \in B_1$,
then by (\ref{EQ2360}) we have $x \rightarrow y$.
Since for $s\in a(x)^\smallfrown a(y)$ we have $y \rightarrow s \rightarrow x$,
by (\ref{EQ2364}) there is $z\in D$ such that $y \rightarrow z \rightarrow x$.
$z\in B_1$
would imply that $y \r z \r x$ and, hence, $y \r x$, which is false.
Thus $z\in A_1$
and, by (\ref{EQ2360}), $x \r z \r y$.

If $x \in A_1$ and $y\in B_1$,
then by (\ref{EQ2360}) we have $y \rightarrow x$.
Since for $s\in x^\smallfrown a(y)$ we have $y \rightarrow s $ and $x\rightarrow s$,
by (\ref{EQ2364}) there is $z\in D$ such that $y \rightarrow z $ and $x\rightarrow z$.
Assuming that $z\in B_1$
we would have $y \r z \r x$
and, hence, $y \r x$, which is false.
Thus $z\in A_1$
and, by (\ref{EQ2360}), $x \r z \r y$.

If $x \in B_1$ and $y\in A_1$,
then by (\ref{EQ2360}) we have $y \rightarrow x$
Since for $s\in a(x)^\smallfrown y$ we have $s \rightarrow x $ and $s\rightarrow y$,
by (\ref{EQ2364}) there is $z\in D$ such that $z \rightarrow x $ and $z\rightarrow y$.
Assuming that $z\in B_1$
we would have $y \r z \r x$, and, hence, $y \r x$, which is false.
Thus $z\in A_1$
and, by (\ref{EQ2360}), $x \r z \r y$.

So $A_1$ is a dense subset of $\la D,\r \upharpoonright D \ra$ and the proof for $B_1$ is similar. This implies that $\la D,\r \upharpoonright D \ra$
is a dense linear order.

Suppose that there are $x=\min _{ \la D,\r \upharpoonright D \ra} D$ and $y=\max_{\la D,\r \upharpoonright D \ra} D$. Then
\begin{equation}\label{EQ2287}
\forall z\in D\setminus \{ x, y\} \;\; x\r z \r y.
\end{equation}
If $x\rightarrow y$, then, since $x\r y$, by (\ref{EQ2360}) we have $x,y\in A$ or $x,y\in B$.
But for $s\in a(y)^\smallfrown x$ we have $s\rightarrow x$ and $s\rightarrow y$ and,
by (\ref{EQ2364}), there is $z\in D$ such that $z\rightarrow x$ and $z\rightarrow y$
which by (\ref{EQ2287})
implies that $x$ and $y$ are in different elements of the partition $\{ A,B\}$ and we have a contradiction.

If $y\rightarrow x$, then, by (\ref{EQ2360}), $x$ and $y$ are in different elements of the partition $\{ A,B\}$.
But for $s\in y^\smallfrown x$ we have $y \rightarrow s \rightarrow x$ and,
by (\ref{EQ2364}), there is $z\in D$ such that $y \rightarrow z\rightarrow x$.
So, by (\ref{EQ2287}) and (\ref{EQ2360}) we have $x,y\in A$ or $x,y\in B$
and we have a contradiction again.
\kdok
Now we prove that $\P (\Y )$ is a dense suborder of  $\P (\S (2))$.
If $D\in \P (\S (2))$,  then, by Claim \ref{T2652}, $\la D ,\r \upharpoonright D \ra$ is a dense linear order
and $A_1:=D\cap A$ and $B_1:=D\cap B$ are its dense subsets.
Let $D'$ be the set obtained from $D$ by deleting its end point, if it exists. Then $\la D' ,\rho \upharpoonright D' \ra $
is a dense linear order without end points, $A_1':= D' \cap A$ and $B_1':= D' \cap B$ are its dense and disjoint subsets and, hence
$\BD ':= \la D' ,\rho \upharpoonright D', A_1',B_1' \ra \models \T _2$, which, since the theory $\T _2$ is $\o$-categorical, implies that
$\BD '\cong \Y$; so $D'\in \P (\Y )$ and, clearly, $D'\subset D$.
Thus $\P (\Y )$ is dense in $\P (\S (2))$ and, hence, $\P (\S (2))\equiv_{forc} \P (\Y )\cong \P (\Q _2)$ so $\P (\S (2))\equiv_{forc}\P (\Q _2)$.

The second statement follows from the first,  Theorem \ref{T4120}(b) and (\ref{EQ2372}).
\hfill $\Box$
\section{The digraphs $\S (3)$, $\BT [\BI_n]$ and $\BI_n [\BT ]$}\label{S4}
\paragraph{The digraph $\S (3)$}
Again we consider the subset $S=: \{ e^{qi}: q\in Q \}$
of the unit circle $S^1$ in the complex plane.
 If  $r: S^1\rightarrow S^1$ is the rotation given by $r(e^{ti})=e^{(t+\frac{2\pi}{3})i}$ and $x=e^{qi}\in S$,
 then $r(x),r^2(x)\not \in S$, where $r^2(x):=r(r(x))$, and the points $x$, $r(x)$ and $r^2 (x)$
are vertices of a equilateral triangle.
It is clear that the $L_b$-structure $\S (3):=\la S, \rightarrow \ra$, where $\rightarrow$ is the binary relation on $S$ defined by
\begin{equation}\label{EQ2368}\textstyle
e^{q_1 i}\rightarrow e^{q_2 i} \;\;\Leftrightarrow \;\; q_2 -q_1 \in \bigcup _{k\in \Z }(2k\pi ,2k\pi +\frac{2\pi}{3}  ),
\end{equation}
is a digraph;
in fact we have $x\rightarrow y$ iff $y\in x^\smallfrown r(x)$, where
for non-antipodal points $s,t\in S ^1$ by $s^\smallfrown t$ we denote the set of elements of $S$
belonging to the shorter arc of $S^1$ determined by $s$ and $t$.
The digraph $\S (3)$ is not a tournament; namely the $L_b$-formula
$\theta (u,v)  :=  u\neq v \land \neg R(u,v) \land \neg R(v,u)$ defines the incomparability
relation, $\parallel$, in $\S (3)$: for $x,y\in S$,
$$
x\parallel y \Leftrightarrow x\neq y \land \neg x\rightarrow y \land \neg y\rightarrow x
$$
and we have $x\parallel y$ iff $y\in r(x)^\smallfrown r^2(x)$. In addition, $y\rightarrow x$ iff $y\in r^2(x)^\smallfrown x$
and, hence, $\{ \Delta _S ,\rightarrow, \rightarrow ^{-1}, \parallel\}$ is a partition of the set $S^2$, where $\Delta _S=\{ \la x,x\ra :x\in S\}$
is the diagonal of $S$. $\S (3)$ is one of continuum may ultrahomogeneous digraphs \cite{Che2}.

For convenience, let $L_3=\la R, \a ,\b ,\g \ra$, where $\ar (R)=2$ and $\ar (\a )=\ar (\b )=\ar (\g )=1$,
and let $L_b := \la R\ra$.
It is evident that  $\{ A , B ,C \}$ is a partition of the set $S$, where
\begin{eqnarray*}\textstyle
A & := &\textstyle\Big\{ e^{qi}: q\in \bigcup _{k\in \Z}(\frac{3\pi}{6}+2k\pi , \frac{7\pi}{6}+2k\pi) \cap Q\Big\} ,\\
B & := &\textstyle\Big\{ e^{qi}: q\in \bigcup _{k\in \Z}(\frac{7\pi}{6}+2k\pi , \frac{11\pi}{6}+2k\pi)\cap Q\Big\} ,\\
C & := &\textstyle\Big\{ e^{qi}: q\in \bigcup _{k\in \Z}(\frac{11\pi}{6}+2k\pi , \frac{15\pi}{6}+2k\pi)\cap Q\Big\} ,
\end{eqnarray*}
and, clearly,
\begin{equation}\label{EQ2405}
\la A , \rightarrow\upharpoonright A \ra \cong \la B , \rightarrow\upharpoonright B \ra \cong \la C , \rightarrow\upharpoonright C \ra \cong \Q ,
\end{equation}
\begin{equation}\label{EQ2402}
\Big((A\times C)\cup (C\times B)\cup (B\times A)\Big)\;\; \cap \;\rightarrow \;\; = \;\emptyset \;\;\mbox{ and}
\end{equation}
\begin{equation}\label{EQ2404}
\Big((C\times A)\cup (B\times C)\cup (A\times B)\Big) \;\;\cap  \;\, \rightarrow ^{-1} \;\;=\; \emptyset .
\end{equation}
Now, $\la S, \rightarrow , A ,B ,C\ra$ is an $L_3$-structure,
the $L_3$-formula
\begin{eqnarray*}
\l (u, v) & := & \Big[ \Big(( \a (u) \land \a (v) ) \lor ( \b (u) \land \b (v))\lor ( \g (u) \land \g (v))\Big)\land R(u,v) \Big]  \nonumber \\
          &\lor& \Big[ \Big(( \a (u) \land \g (v))\lor (\g (u) \land \b (v))\lor (\b (u) \land \a (v))\Big)\land R(v,u) \Big]\\
          &\lor& \Big[ \Big(( \g (u) \land \a (v))\lor (\b (u) \land \g (v))\lor (\a (u) \land \b (v))\Big)\land \;\theta (u,v) \Big]
\end{eqnarray*}
defines a new binary relation $\tau$ on $S$
\begin{eqnarray}
\tau & = &\Big[\Big((A\times A)\cup (B\times B)\cup (C\times C)\Big) \;\;\cap \;  \rightarrow   \Big] \nonumber \\
     & \cup & \Big[\Big((A\times C)\cup (C\times B)\cup (B\times A)\Big)\;\; \cap \;\rightarrow ^{-1} \Big] \label{EQ2401} \\
     & \cup & \Big[\Big((C\times A)\cup (B\times C)\cup (A\times B)\Big) \;\;\cap  \;\, \parallel   \Big] \nonumber
\end{eqnarray}
and $\la S, \tau , A, B, C \ra$ is an $L_3$-structure as well. By (\ref{EQ2401}) we have
\begin{eqnarray}
\tau ^{-1}& = &\Big[\Big((A\times A)\cup (B\times B)\cup (C\times C)\Big) \;\;\cap \;  \rightarrow ^{-1}  \Big] \nonumber \\
          & \cup & \Big[\Big((C\times A)\cup (B\times C)\cup (A\times B)\Big) \;\;\cap  \;\, \rightarrow   \Big]\label{EQ2403}\\
     & \cup & \Big[\Big((A\times C)\cup (C\times B)\cup (B\times A)\Big)\;\; \cap \;\parallel \Big]  \nonumber.
\end{eqnarray}
For completeness we include a proof of the following well-known fact.
\begin{fac}\label{T2401}

(a) $\la S, \tau , A, B, C \ra \cong \Q _{3}$;

(b)  $\la S, \rightarrow , A, B, C \ra$ and $\la S, \tau , A, B, C \ra$ are $\Sigma _0$-bi-definable $L_3$-structures.
\end{fac}
\dok
(a) Since $\rightarrow$, $\rightarrow ^{-1}$ and $\parallel$ are irreflexive and pairwise disjoint binary relations on $S$,
by (\ref{EQ2401})  the relation $\t$ is irreflexive and,
by (\ref{EQ2401}) and (\ref{EQ2403}), $\tau \cap \tau ^{-1}=\emptyset$; so the relation $\t$ is asymmetric; so $\la S,\tau \ra$ is a digraph.
In addition, by (\ref{EQ2405}) - (\ref{EQ2403})  we have $\tau \cup \tau ^{-1}=S^2 \setminus \Delta _S $,
which means that $\la S,\tau \ra$ is a tournament.

Suppose that the relation $\tau$ is not transitive. Then $x \t y \t  z\t x$, for some $x,y,z\in S$, and, by (\ref{EQ2405}),
$x$, $y$ and $z$ are not in the same of the sets $A$, $B$ and $C$.

Suppose that two of these points belong to one of these sets, say $x,y\in A$, which implies that $x\rightarrow y$.
If $z\in B$, then, by (\ref{EQ2401}), $y\parallel z$ and $x\rightarrow z$ and, hence $y,z\in x^\smallfrown r(x)$, which implies that $y\not\parallel z$
and we have a contradiction.
If $z\in C$, then, by (\ref{EQ2401}), $z\parallel x$ and $z\rightarrow y$ and, hence $x,z\in r^2(y)^\smallfrown y$, which implies that $x\not\parallel z$
and we have a contradiction.
In a similar way we show that whenever two of the points belong to one of the elements of the partition we obtain a contradiction.

Thus $x$, $y$ and $z$ are in different elements of the partition  and by (\ref{EQ2401}) we have:
if $\la x,y,z\ra\in (A\times C\times B)\cup (B\times A\times C)\cup (C\times B\times A)$, then
$x\rightarrow z\rightarrow y\rightarrow x $ so $\{ x,y,z\}$ is a copy of the oriented triangle, $C_3$,  in $\S (3)$, which is impossible;
if $\la x,y,z\ra\in (A\times B\times C)\cup (B\times C\times A)\cup (C\times A\times B)$, then
$x\parallel z\parallel y\parallel x $ and $\{ x,y,z\}$ is a copy of the empty digraph, $E_3$, in $\S (3)$, which is impossible again.

A proof that $A$, $B$ and $C$ are dense sets in the linear order $\la S,\tau\ra$ follows from the proof of Claim  \ref{T2655} (take $D=S$).
Suppose that $m=\min S$ and, say $m\in A$;
but by (\ref{EQ2401}) and (\ref{EQ2405}) we have $\la A , \tau \upharpoonright A \ra =\la A , \rightarrow\upharpoonright A \ra \cong \Q$ and this is impossible. So
$\la S,\tau\ra$ is a dense linear order without end points, $\la S, \tau , A, B, C \ra \models \T_3$ and, hence, $\la S, \tau , A, B, C \ra \cong \Q _{3}$.

(b) First, $\tau =\{ \la x,y \ra\in S^2 :\la S, \rightarrow , A, B, C \ra \models \l [x,y] \}$ and we show that
$\rightarrow  \;=\{ \la x,y \ra\in S^2 :\la S, \tau , A, B, C \ra \models \mu [x,y] \}$, where $\mu (u,v)$ is the $L_3$-formula
\begin{eqnarray*}
\mu (u, v) & := & \Big[ \Big(( \a (u) \land \a (v) ) \lor ( \b (u) \land \b (v))\lor ( \g (u) \land \g (v))\Big)\land \;R(u,v) \Big]  \nonumber \\
           &\lor& \Big[ \Big(( \g (u) \land \a (v))\lor (\b (u) \land \g (v))\lor (\a (u) \land \b (v))\Big)\land \neg R(u,v) \Big],
\end{eqnarray*}
that is, defining $U:= A^2\cup B^2\cup C^2$, $V:=(C\times A)\cup (B\times C)\cup (A\times B)$ and $W:=(A\times C)\cup (C\times B)\cup (B\times A)$ we prove that
\begin{equation}\label{EQ2406}
\rightarrow = (U \cap \tau ) \cup  (V \setminus \tau   ) .
\end{equation}
By (\ref{EQ2402}) we have $\rightarrow = (U \cap \rightarrow ) \cup  (V \cap \rightarrow )$ and, by (\ref{EQ2401}), $U \cap \tau =U \cap \rightarrow $.
By (\ref{EQ2401}) and (\ref{EQ2404})  we have
$V \setminus \tau =V \setminus \parallel=V \cap (\rightarrow \cup \rightarrow ^{-1})=V \cap \rightarrow$ so (\ref{EQ2406}) is true.
Since the formulas $\l$ and $\mu$ are quantifier free, statement (b) is proved.
\hfill $\Box$
\begin{te}\label{T2657}
$\P (\Q _3)$  densely embeds in $\P (\S (3))$ and, hence, $\B _{\S (3)}\cong \B _{\Q _3}$.

Each countable digraph equimorphic with $\S (3)$ has property $\CP _1$.
\end{te}
\dok
Let $\Y :=\la S, \tau , A ,B, C \ra$. By Fact \ref{T2401}(a) we have $\P (\Q _3) \cong \P (\Y)$ so
it is sufficient to show that $\P (\Y )$ is a dense subset of $\P (\S (3))$.
We prove first that $\P (\Y ) \subset \P (\S (3))$.
So, if $D\in \P (\Y )$, then there is an isomorphism
\begin{equation}\label{EQ2362}
F: \la S, \tau , A ,B, C \ra \rightarrow_{iso} \la D, \tau \upharpoonright D , A\cap D ,B\cap D , C\cap D \ra
\end{equation}
and in order to prove that $D\in \P (\S (3))$ it remains to be shown that the mapping
$F: \la S, \rightarrow \ra \rightarrow \la D, \rightarrow \upharpoonright D \ra$ is an isomorphism.
By Fact \ref{T2401}(b), the relation $\rightarrow$ is defined by the $L_3$-formula $\mu$ in the structure $\Y$, that is
\begin{equation}\label{EQ2361}
\forall x,y  \in S \;\; \Big(x \rightarrow y \;\Leftrightarrow \;\la S, \tau , A,B,C \ra \models \mu [x,y] \Big).
\end{equation}
Now for $x,y \in S$ we have:
$x\rightarrow y$
iff (by (\ref{EQ2361})) $\la S, \tau , A ,B ,C \ra  \models \mu [x,y] $
iff (by (\ref{EQ2362})) $\la D, \tau \upharpoonright D , A\cap D ,B\cap D , C\cap D \ra  \models \mu [F(x),F(y)] $
iff (since $\mu$ is a $\Sigma _0$-formula and, thus, $(D,S)$-absolute) $\la S, \tau , A ,B ,C \ra   \models \mu [F(x),F(y)] $
iff (by (\ref{EQ2361})) $F(x)\rightarrow F(y)$.
Thus $F: \la S, \rightarrow \ra \rightarrow \la D, \rightarrow \upharpoonright D \ra$ is an isomorphism,
$D\in \P (\S(3))$ and  $\P (\Y ) \subset \P (\S (3))$ indeed.
\begin{cla}\label{T2655}
If $D\in \P (\S (3))$, then $\la D ,\tau \upharpoonright D\ra$ is a dense linear order and the sets
$A_1:=A\cap D$, $B_1:=B\cap D$ and $C_1:=C\cap D$
are dense in $\la D ,\tau \upharpoonright D\ra$.
\end{cla}
\dok
By Fact \ref{T2658}(b), if $D\in \P (\S (3))$, then $\BD:= \la D,\rightarrow \upharpoonright D\ra \prec\S(3)$. So, by the Tarski-Vaught theorem,
for each $L_b$-formula $\theta (u,v,w)$ we have:
\begin{equation}\label{EQ2363}
\forall x,y \in D \;
\Big( \exists s\in S \;\; \S(3) \models  \theta [x,y,s]
\;\Rightarrow \exists z\in D \;\;
\BD \models \theta [x,y,z]   \Big).
\end{equation}
By Fact \ref{T2401}(a) $\la D,\tau \upharpoonright D\ra$ is a linear order and we prove that $A_1$ is its dense subset.
So, assuming that $x,y \in D$ and $x\t y$ we will find a $z\in A_1$ such that $x\t  z\t y$.

If $x,y \in A_1$,
then by (\ref{EQ2401}) we have $x\rightarrow y$.
Since for $s\in x^\smallfrown y$ we have $x\rightarrow s \rightarrow y$, by (\ref{EQ2363}), there is $z\in D$ such that
$x \rightarrow z \rightarrow y$.
Since $x,y \in A$, by (\ref{EQ2402}) we have $z\not\in B\cup C$,
which implies that $z\in A_1$.
Thus, by (\ref{EQ2401}) we have $x\t  z\t y$.

If $x,y \in B_1$,
then by (\ref{EQ2401}) we have $x\rightarrow y$. Since for $s\in r^2(x)^\smallfrown r^2(y)$ we have $s \rightarrow x$ and $y\parallel s$,
by (\ref{EQ2363}) there is $z\in D$ such that $z \rightarrow x$ and $y\parallel z$.
Since $x \in B$, by (\ref{EQ2402}) we have $z\not\in  C$, and assuming that $z\in B$ we would have $y\not\parallel z$
(because $\la B,\rightarrow \upharpoonright B\ra$ is a linear order).
Thus $z\in A_1$ and, by (\ref{EQ2401}), $\la x,z\ra\in (B\times A)\cap \rightarrow ^{-1}\subset \t$ and $\la z,y\ra\in (A\times B)\cap \parallel \subset \t$.
So we have $x\t  z\t y$.

If $x,y \in C_1$, then by (\ref{EQ2401}) we have $x\rightarrow y$. Since for $s\in r(x)^\smallfrown r(y)$ we have $x\parallel s$  and $y \rightarrow s$,
by (\ref{EQ2363}) there is $z\in D$ such that $x\parallel z$  and $y \rightarrow z$.
Since $y \in C$, by (\ref{EQ2402}) we have $z\not\in  B$,
and assuming that $z\in C$ we would have $x\not\parallel z$
(because $\la C,\rightarrow \upharpoonright C\ra$ is a linear order).
Thus $z\in A_1$ and, by (\ref{EQ2401}), $\la x,z\ra\in (C\times A)\cap \parallel \subset \t$ and $\la z,y \ra\in (A\times C)\cap \rightarrow ^{-1}\subset \t$.
So we have $x\t  z\t y$.

If $x \in A_1$, $y \in B_1$, then by (\ref{EQ2401}) we have $x\parallel y$. Since for $s\in x^\smallfrown r^2(y)$ we have  $x \rightarrow s $ and $s\parallel y$,
by (\ref{EQ2363}) there is $z\in D$ such that  $x \rightarrow z $ and $z\parallel y$.
Since $x \in A$, by (\ref{EQ2402}) we have $z\not\in C$,
and assuming that $z\in B$ we would have $z\not\parallel y$
(because $\la B,\rightarrow \upharpoonright B\ra$ is a linear order).
Thus $z\in A_1$ and, by (\ref{EQ2401}),
$\la x,z\ra\in (A\times A)\cap \rightarrow \subset \t$ and
$\la z,y \ra\in (A\times B)\cap \parallel \subset \t$.
So we have $x\t  z\t y$.

If $x \in A_1$, $y \in C_1$, then by (\ref{EQ2401}) we have $y\rightarrow x$. Since for $s\in x^\smallfrown r(y)$ we have $x \rightarrow s $ and $y \rightarrow s $,
by (\ref{EQ2363}) there is $z\in D$ such that $x \rightarrow z $ and $y \rightarrow z $.
Since $x \in A$, by (\ref{EQ2402}) we have $z\not\in C$; since $y \in C$, by (\ref{EQ2402}) we have $z\not\in B$,
Thus $z\in A_1$ and, by (\ref{EQ2401}),
$\la x,z\ra\in (A\times A)\cap \rightarrow \subset \t$ and
$\la z,y \ra\in (A\times C)\cap \rightarrow ^{-1} \subset \t$.
So, $x\t  z\t y$.

If $x \in B_1$, $y \in C_1$, then by (\ref{EQ2401}) we have $x\parallel y$. Since for $s\in r^2(x)^\smallfrown r(y)$ we have $s \rightarrow x$ and $y \rightarrow s$,
by (\ref{EQ2363}) there is $z\in D$ such that $z \rightarrow x $ and $y \rightarrow z $.
Since $x \in B$, by (\ref{EQ2402}) we have $z\not\in C$; since $y \in C$, by (\ref{EQ2402}) we have $z\not\in B$,
Thus $z\in A_1$ and, by (\ref{EQ2401}),
$\la x,z\ra\in (B\times A)\cap \rightarrow  ^{-1}\subset \t$ and
$\la z,y \ra\in (A\times C)\cap \rightarrow ^{-1} \subset \t$.
Thus, $x\t  z\t y$.

If $x \in B_1$, $y \in A_1$, then by (\ref{EQ2401}) we have $y\rightarrow x$. Since for $s\in r^2(x)^\smallfrown y$ we have $s \rightarrow x $ and $s \rightarrow y$,
by (\ref{EQ2363}) there is $z\in D$ such that $z \rightarrow x $ and $z \rightarrow y $.
Since $x \in B$, by (\ref{EQ2402}) we have $z\not\in C$,
and since $y \in A$, by (\ref{EQ2402}) we have $z\not\in B$.
Thus $z\in A_1$ and, by (\ref{EQ2401}),
$\la x,z\ra\in (B\times A)\cap \rightarrow^{-1} \subset \t$ and
$\la z,y \ra\in (A\times A)\cap \rightarrow \subset \t$.
So, $x\t  z\t y$.

If $x \in C_1$, $y \in A_1$, then by (\ref{EQ2401}) we have $x\parallel y$. Since for $s\in r(x)^\smallfrown y$ we have $x \parallel s $ and $s \rightarrow y $,
by (\ref{EQ2363}) there is $z\in D$ such that $x \parallel z $ and $z \rightarrow y $.
Since $y \in A$, by (\ref{EQ2402}) we have $z\not\in B$;
and assuming that $z\in C$ we would have $x\not\parallel z$
(because $\la C,\rightarrow \upharpoonright C\ra$ is a linear order).
Thus $z\in A_1$ and, by (\ref{EQ2401}),
$\la x,z\ra\in (C\times A)\cap \parallel \subset \t$ and
$\la z,y \ra\in (A\times A)\cap \rightarrow  \subset \t$.
So we have $x\t  z\t y$.

If $x \in C_1$, $y \in B_1$, then by (\ref{EQ2401}) we have $y\rightarrow x$. Since for $s\in r(x)^\smallfrown r^2(y)$ we have $x \parallel s$ and $y \parallel s$,
by (\ref{EQ2363}) there is $z\in D$ such that $x \parallel z $ and $y \parallel z$.
Since $\la C,\rightarrow \upharpoonright C\ra$ and $\la B,\rightarrow \upharpoonright B\ra$ are linear orders,
assuming that $z\in C$ (resp.\ $z\in B$) we would have $x\not\parallel z$ (resp.\ $y\not\parallel z$).
Thus $z\in A_1$ and, by (\ref{EQ2401}),
$\la x,z\ra\in (C\times A)\cap \parallel\subset \t$ and
$\la z,y \ra\in (A\times B)\cap \parallel \subset \t$.
So we have $x\t  z\t y$.

Proofs that $B_1$ and $C_1$ are dense sets in the linear order $\la D,\tau \ra$ are similar.
\kdok

Now, if $D\in \P (\S (3))$, then, by Claim \ref{T2655}, $\la D ,\tau \upharpoonright D \ra$ is a dense linear order
and $A\cap D$, $B\cap D$ and  $C\cap D$ are dense sets in $\la D ,\tau \ra$.
Let $D'$ be the set obtained from $D$ by deleting its end points, if they exist.
Then $\la D' ,\tau \upharpoonright D' \ra $ is a dense linear order without end points
and $\{ A\cap D', B\cap D',C\cap D'\}$ is a partition of $D'$ into three dense subsets of $\la D' ,\tau \upharpoonright D' \ra $.
Thus ${\mathbb D }':=\la D' ,\tau \upharpoonright D', A\cap D', B\cap D',C\cap D'\ra$ is a substructure of $\Y$
and ${\mathbb D }' \models \T _3$,
which, since the theory $\T _3$ is $\o$-categorical and, by Fact \ref{T2401}(a), $\Y \models \T _3$,
implies that ${\mathbb D }'\cong \Y$.
So $D'\in \P (\Y )$, $D'\subset D$ and $\P (\Y )$ is a dense suborder of  $\P (\S (3))$ indeed.
Thus $\P (\S (3))\equiv_{forc} \P (\Y )\cong \P (\Q _3)$ and, hence, $\P (\S (3))\equiv_{forc}\P (\Q _3)$.

The second statement follows from the first,  Theorem \ref{T4120}(b) and (\ref{EQ2372}).
\hfill $\Box$
\paragraph{Wreath products $\BT [\BI_n]$ and $\BI_n [\BT ]$.}
One subclass of the class of all ultrahomogeneous digraphs (Cherlin's list \cite{Che2}) is described as follows.
Let $\BT$ be an ultrahomogeneous tournament (thus $\BT \in \{ \Q , \BT ^\infty, \S (2)\}$)
and, for an integer $n\geq 2$, let $\BI_n$ denote the digraph with $n$ vertices and with no arrows.
Then the digraphs

- $\BT [\BI_n]$ (obtained by replacement of each point of $\BT$ by a copy of $\BI _n$) and

- $\BI_n [\BT ]$ (obtained by replacement of each point of $\BI _n$ by a copy of $\BT$)

\noindent
are ultrahomogeneous and imprimitive.
More precisely, the $L_b$-formula $\f (u,v):= \neg R(u,v)\land \neg R(v,u)$ defines the ``unrelatedness" binary relation $\sim$ on the domain,
and, hence, all automorphisms preserve it.

It is easy to see that all embeddings of $\BT [\BI_n]=\bigcup _{t\in T}I_n^t$ preserve the relation $\sim$ as well and hence,
$\P (\BT [\BI_n])=\{ \bigcup _{t\in A }I_n^t: A\in \P (\BT )\} \cong \P (\BT )$. So, the digraphs $\Q [\BI_n]$ and $\S (2)[\BI_n]$
have property $\CP _1$ while $ \BT ^\infty [\BI_n]$ has $\CP _2$.

On the other hand, the digraphs $\BI_n [\BT ]$ are disconnected and, by Theorem 5.2 of \cite{Kstr},
$\P (\BI_n [\BT ])\cong \P (\BT )^n$. Thus, for example, the poset $\P (\BI_n [\S (2)  ])\equiv _{forc}(\S \ast \pi )^n$.

\vspace{2mm}
\noindent
{\bf Acknowledgement.}
Both authors are supported by the Science Fund of the
Republic of Serbia, Program IDEAS, Grant No. 7750027: Set-theoretic, model-theoretic
and Ramsey-theoretic phenomena in mathematical structures: similarity
and diversity--SMART.
In addition, the second author is partially
supported by grants from NSERC (455916) and CNRS (IMJ-PRG UMR7586).

\end{document}